\documentclass[12pt]{article}
\usepackage[T2A]{fontenc}
\usepackage[cp866]{inputenc}
\usepackage[english,russian]{babel}
\usepackage[tbtags]{amsmath}
\usepackage{amsfonts,amssymb,mathrsfs,amscd,amsthm,hyper}
\theoremstyle{main}
\newtheorem{theorem}{Theorem}
\newtheorem{theoremSta}{Theorem of Stahl}
\newtheorem{theoremGoRa}{Theorem of Gonchar--Rakhmanov}
\newtheorem{theoremGol}{Theorem of Goluzin}
\newtheorem{theoremGoRacon}{Corollary of Gonchar--Rakhmanov theorem}
\newtheorem{theoremBaStYa}{Theorem of Baratchart--Stahl-Yattselev}

\def\Re{\operatorname{Re}}
\def\Im{\operatorname{Im}}

\def\const{\operatorname{const}}

\def\mcap {\operatorname{cap}}

\def\supp{\operatorname{supp}}

\def\equ{\operatorname{eq}}
\def\({\left(}
\def\){\right)}
\def\HH{\mathcal H}
\def\myVV{\mathfrak V}

\def\myA{\mathcal A}
\def\myk{\mathfrak K}
\def\MM{M}
\def\CC{\mathbb C}

\let\<\langle
\let\>\rangle

\let\myo\overline
\let\myt\widetilde

\def\bad{\spaceskip=0.33emplus0.6emminus0.15em\immediate\write5{\string\bad}}

\let\geq\geqslant
\let\leq\leqslant

\def\eqref#1{(\ref{#1})}

\theoremstyle{remark}

\begin{document}


\title{On the existence of compacta of minimal capacity in the theory of
rational approximation of multi-valued analytic functions}


\date{22.05.2015}

\author{Viktor~I.~Buslaev, Sergey~P.~Suetin}

\maketitle

\begin{abstract}

For an interval $E=[a,b]$ on the real line, let $\mu$ be either the equilibrium measure, or the normalized Lebesgue measure of $E$, and let
$V^{\mu}$ denote the associated logarithmic potential. In the present paper, we construct a
function~$f$ which is analytic on $E$ and possesses four branch points of
second order outside of $E$ such that the family of the admissible compacta of
$f$ has no minimizing elements with regard to the extremal theoretic-potential
problem, in the external field equals $V^{-\mu}$.

Bibliography: 35~items.
\end{abstract}

rational approximants, Pad\'{e} approximation, orthogonal
polynomials, distribution of poles, convergence in capacity

\footnotetext{This work is supported by the Russian Science Foundation (RSF)
under the grant 14-21-00025.}



\begin{center}\it
Dedicated to the memory of Andrei Aleksandrovich Gonchar and Herbert
Stahl
\end{center}


\section{Notations}\label{s1}

Throughout the paper, we use the following notations.

$\MM (K)$ -- the space of all positive unit Borel measures $\mu$ with supports
$S(\mu)=\supp\mu\subset K$, where $K$ is a compact set, $K\subset\overline{\Bbb C}$.

$\delta_z$ -- the Dirac measure at a point $z$.

For a finite set $E=\{e_1,\dots ,e_n\}$ of points (counted with their
multiplicities) on $\myo{\Bbb C}$, we introduce the measure
$\delta_E:=\sum_{k=1}^n\delta_{e_k}$.

For a polynomial $Q$ of degree $n$, we denote
$\delta_Q:=\delta_{\{ z_{1},\dots ,z_{n}\}}$, where $\{ z_{1},\dots
,z_{n}\}$ are the zeros of $Q$.

We define a spherically normalized potential ${\mathcal V}^{\mu }$ of the
measure $\mu$ by
\begin{equation}
{\mathcal V}^{\mu }(z)=\int_{|t|\leq 1}\log\frac 1{|z-t|}d\mu (t)+\int_{|t|> 1}\log\frac 1{|1-z/t|}d\mu (t).
\label{pt173}
\end{equation}

$\mcap_\psi K$ denotes the capacity of the compactum $K$ in the presence of
the harmonic external field $\psi$ (the so-called $\psi$-weighted capacity).
It is known ~\cite{SaTo97} that
$\mcap_\psi K$ coincides with the transfinite diameter $d_\psi K$ of
$K$ in the field $\psi$, that is
\begin{equation}
d_\psi K:=\lim_{n\to\infty}
\biggl (\max_{z_1,\dots ,z_n\in K}\prod_{1\leq q<r\leq n}|z_q-z_r|e^{-\big (\psi (z_q)+\psi (z_r)\big)}\biggl)^{\frac 2{(n-1)n}}.
\label{pt172}
\end{equation}

Let $\psi ={\mathcal V}^{-\mu}$, $S(\mu)\cap K=\varnothing$.
For the sake of simplicity we write
$\mcap_\mu K:=\mcap_{{\mathcal V}^{-\mu}}K$.

Consider a sequence of compacta $\{K_n\},\;n=1,2,\dots$. It is
known~\cite{Rak12} that if each $K_n$
is a union of a finite number of continua (the number does not depend on the index $n$),
such that $\lim_{n\to\infty}K_n=K$
(in the Hausdorff metric, see~\eqref{K1}--\eqref{K} below) and if
$S(\mu)\cap K=\varnothing$, then
\begin{equation}
\lim_{n\to\infty}\mcap_\mu K_n=\mcap_\mu K.
\label{h3}
\end{equation}
$\mcap K$ denotes the standard capacity of the compactum $K\subset\Bbb C$ in
the absence of an external field, i.e. when $\psi\equiv 0$. Since ${\mathcal
V}^{\delta_\infty}(z) =0$ for all $z\in\Bbb C$, then $\mcap K=\mcap_{\delta_\infty}K$.

Given an open set $\Omega$, we say that
a sequence of functions $\{R_n\}_{n=1}^\infty$ converges in capacity to the
function $f$ on compact subsets of $\Omega$, if for every $ K\subset\Omega$ and
each $\varepsilon >0$
\begin{equation}
\lim_{n\to\infty}\mcap\{ z\in K:\,\big |(R_n-f)(z)\big |>\varepsilon\} =0;
\label{h2}
\end{equation}
we use
the notation
$$
R_n\overset{\mcap}\rightarrow f,\quad z\in\Omega,\quad n\to\infty.
$$
The notation $\mu_n\overset{_*}\to\mu$
is used for weak convergence of
a sequence of measures $\mu_n$ to the measure $\mu$ as $n\to\infty$.
Recall that
\begin{equation}
\mu_n\overset{_*}\to\mu\quad\text{as}\quad n\to\infty\,\Longleftrightarrow\,\lim_{n\to\infty}\int fd\mu_n=\int fd\mu,
\label{h4}
\end{equation}
for each function $f$ continuous
on $\overline{\Bbb C}$. Also recall that
from each sequence of unit measures one can extract a convergent subsequence.

We denote by $g_K(z,\zeta)$ Green's function of the complement $\overline{\Bbb C}\setminus K$
of the compactum $K$, with a pole at $\zeta\in\overline{\Bbb C}\setminus K.$
(We let $g_K(z,\zeta)=0$, if the points $z$ and $\zeta$ are contained in
two
separate components of $\overline{\Bbb C}\setminus K$.)

$G_K^\mu (z)=\int g_K(z,t)d\mu (t)$ is Green's potential of the measure
$\mu$ whose support $S(\mu)$ does not intersect the compactum $K$.

Now let $S(\mu)\cap K=\varnothing$. Let us remind the notion of the balayage
$\widetilde{\mu }_K$ of a measure $\mu$ onto the compactum $K$. By definition,
the balayage is the
unique measure in the space $\MM (K)$ for which the equality
\begin{equation}
G_K^\mu (z)={\mathcal V}^{\mu -{\widetilde{\mu }_K}}(z)+w^{\mu }_K,
\quad z\in\overline{\Bbb C},
\label{pt77}
\end{equation}
holds quasi everywhere, i.e. except for a set of zero inner
capacity;
$w_K^{\mu }$ is a constant. In the particular case when $\mu =\delta_\infty$,
from equality~\eqref{pt77} we obtain that
\begin{equation}
g_K(z,\infty)={\mathcal V}^{-\lambda_K}(z)+w_K,\quad z\in
\overline{\Bbb C},
\label{pt78}
\end{equation}
where $\lambda_K$ is the equilibrium measure of the compact set $K$ (the balayage of the measure
$\mu =\delta_\infty$ onto $K$) and $w_K$ is the Robin constant for $K.$

We denote ${\mathcal F}$ the class of compacta $F$ of the form
$F=\bigcup_{s=1}^p\overline{\gamma}_s$, where each $\gamma_s$ is
an open analytic curve such that the set
$\overline{\Bbb C}\setminus\overline{\gamma}_s$ is connected and $\gamma_s\cap\overline{\gamma }_t=\varnothing$ for every $s,t=1,\dots ,p$,
$s\neq t.$ If $F=\bigcup_{s=1}^p\overline{\gamma}_s\in {\mathcal F}$, then we use the notation
$F_0$ for the set $\bigcup_{s=1}^p\gamma_s$.

Assume that $S(\mu)\cap F=\varnothing$. We say that
the compactum $F$ is {\it symmetric in the external field} ${\mathcal
V}^{-\mu}$, if $F\in {\mathcal F}$
and the equality
\begin{equation}
\frac{\partial G_F^{\mu}(z)}{\partial n_+}=\frac{\partial G_F^{\mu}(z)}{\partial n_-},
\quad z\in F_0,
\label{S2}
\end{equation} holds, where the normal vectors $n_+$ and $n_-$ are taken in opposite directions on $F$.
In other words, $K$ is said to be an $S$-compactum in this field.

For the purpose of the present paper, there is no need to provide the
general concept of symmetry of a compactum located in an arbitrary harmonic
external field $\psi$ (cf.~\cite{GoRa87}).


${\mathcal E}_m$ stands for the class of compacta $E\subset\overline{\Bbb C}$ of the form
$E=\bigsqcup_{j=1}^mE_j$, where $E_1,\dots ,E_m$ are pairwise disjoint
continua in $\overline{\Bbb C}$ (some of
the latter
may consists of a single point).

$\myA (E)$ represents the class of functions $f$, defined on
$E=\bigsqcup_{j=1}^mE_j$ such that
each restriction $f_j=f|_{E_j}$,
$j=1,\dots ,m$, is holomorphic on $E_j$ and admits an analytic continuation along
every path in $\overline{\Bbb C}$, not passing through a finite point set
$A_{f_j}$, where $A_{f_j}$ contains at least one branch point of $f_j$.

Set $A_f=\bigcup_{j=1}^mA_{f_j}$.

Let $E=\bigsqcup_{j=1}^mE_j\in {\mathcal E}_m$ and $f\in\myA (E)$. We denote
by $\myk_{E,f}$ the family of compacta $K\subset\overline{\Bbb C}$
such that $K\cap E=\varnothing$, and each function $f_j=f|_{E_j}$ is holomorphic
(i.e., analytic and single-valued) in that connected component
$D_j$ of
$\overline{\Bbb C}\setminus K$, which contains $E_j$; $f_j=f_k$ if $D_j=D_k$.
We shall refer to these compacta as admissible compacta for the function
$f\in\myA (E)$.

\section{Statement of the problem and discussion }\label{s2}

Let $E\in {\mathcal E}_m$, $f\in\myA (E)$, and $\mu\in\MM (E)$.
Recently, special attention has been
paid to the existence problem of an admissible compactum $F$
minimizing the weighted capacity; i.e., whether there exists
$F\in\myk_{E,f}$ such that
\begin{equation}
\mcap_{\mu }F=\inf_{K\in\myk_{E,f}}\mcap_{\mu }K .
\label{q1}
\end{equation}

If the problem is solvable, then we say that the class of compacta
$\myk_{E,f}$ contains a
minimizing element
with respect to the measure $\mu$ (or,
equivalently, to the field ${\mathcal V}^{-\mu }$). Otherwise, there is no
minimizing element with respect to $\mu$. If such an element exists, for brevity let
us
call it a
`$\mu$-minimizing element'.
In the
particular case, when
\begin{equation}
E=\{\infty\},\quad\mu=\delta_\infty,\quad
f(z)=\sqrt[q]{(z-a_1)\dotsb(z-a_q)}-z,\quad f\in H(\{\infty\}),
\label{q2}
\end{equation}
the above mentioned problem coincides with the well known Chebotarev
problem on the existence of a continuum of minimal (standard) capacity among all
continua which contain the points $a_1,\dots ,a_q$. Chebotarev's problem was
solved in the 1930's by G.~Gr\"otzsch, and independently by
M.~A.~Lavrentiev. G.~Gr\"otzsch~\cite{Gro30} carried out the proof of the
uniqueness
of the extremal continuum, following his
`strip method'. Under the additional basic assumption that the
continuum in question is a union of analytic curves, M.~A.~Lavrentiev \cite{Lav30},
~\cite{Lav34} solved Chebotarev's problem
applying a variation-geometrical approach.

In a series of articles
from 1946 through 1951,
G.~M.~Goluzin created his own
method, different from the earlier approach introduced by M.~Schiffer;
namely,
the so-called method of internal variations. Applying his method, G. M. Goluzin
solved Chebotarev's problem in terms of quadratic differentials
(see~\cite{Kuz06}).

\begin{theoremGol}[(see~\cite{Kuz06})]
Under conditions (10), the family of compacta $\myk
_{\{\infty\},f}$ possesses a unique
minimizing element
$F$ with respect to the measure $\delta_\infty$ (or, equivalently, to the field
${\mathcal V}^{-\delta_\infty }=0$). The continuum $F$ does not separate the
plane; furthermore, $F$ is a union of the closures of the critical trajectories
of the quadratic differential $\frac{-B(z)}{A(z)}dz^2$, where $A(z)$ and $B(z)$
are monic polynomials, $A(z)=(z-a_1)\dots (z-a_q),$ $B(z)$ is of degree $q-2$ and
uniquely determined by the polynomial $A$.
Green's function $g(z,\infty)$ of the domain $\overline{\Bbb C}\setminus F$ with pole at infinity is given by
$$
g(z,\infty)=\Re\int_a^z\sqrt{B(t)/A(t)}dt.
$$
\end{theoremGol}

The information about
zeros of the polynomial $B(z)$ results from conditions related to the connectivity of the union of the closures of the critical trajectories of the quadratic differentials
$\frac{-B(z)}{A(z)}dz^2$. The
explicit evaluation
of $B(z)$ is a very hard problem studied, so far, only in the case
$q=3$. For details on Chebotarev's problem,
and some other
related problems
in the theory of functions, as well, the reader
may be
referred to~\cite{Kuz80},~\cite{Kuz06},~\cite{Str84}.

Note that Chebotarev's problem concerns a concrete function $f$~\eqref{q2}. In contrast to it, in the
eighties of the last century H.~Stahl gave a positive answer to the general problem of existence of a minimizing element
in the class $\myk_{\{\infty\},f}$ for an arbitrary function $f$ in $\myA (\{\infty\})$.
Furthermore, Stahl
described the extremal compactum $F$ in terms of a quadratic differential.
He established its symmetry property
(in the absence of an external field) and using the symmetry, deduced that
the classical diagonal Pad{\'e} approximants converge in capacity to the
function $f$ in the domain $D:=\overline{\Bbb C}\setminus F$. The latter is known
as Stahl's domain; see~\cite{Sta85a}--\cite{Sta86b}.

We remind the reader of the definition of
the classical diagonal Pad{\'e} approximant $R_n$ of a function $f\in
\HH(\{\infty\})$; that is $R_n=P_{n}/Q_{n}$, where
\begin{equation}
\deg P_{n}\leq n,\quad\deg Q_{n}\leq n,\quad Q_{n}\not\equiv 0,
\label{h1}
\end{equation}
and
\begin{equation}
\label{pade}
(Q_{n}f-P_{n})(z)=O\(\frac 1{z^{n+1}}\),\quad z\to\infty .
\end{equation}

Here, and in
the sequel,
we
denote
$\HH(E)$ the class of functions
holomorphic on $E$ (i.e., in some neighborhood of $E$).

\begin{theoremSta}
Let $f\in\myA (\{\infty\})$. Then there is a compactum
$F\in\myk_{\{\infty\},f}$ such that

$1^\circ$.
$F$ is of minimal capacity; i.e.,
$$
\quad\quad\quad\quad\quad\quad\quad\quad\quad\mcap
F =\min_{K\in\myk_{\{\infty\},f}}\mcap K\,.
\quad\quad\quad\quad\quad\quad\quad\quad\quad\quad\quad\quad $$

$2^\circ$.
$F$ is symmetric (in the field ${\mathcal V}^{-\delta
_\infty}=0$), and the complement $\overline{\Bbb C}\setminus F$ is connected.

$3^\circ$. Let $R_n=P_n/Q_n$, $n=1,2,\dots$, be the classical diagonal Pad{\'e} approximants
of $f$. Then
$$
\frac 1n\delta_{Q_{n}}\overset{_*}\to\mu^{\equ}_F\quad \text{and}\quad
R_n\overset{\mcap}\rightarrow f,\quad z\in\overline{\Bbb C}\setminus F, \hbox{
as }\, n\to\infty,
$$
where $\mu^{\equ}_F$ is the equilibrium measure of $F$.
The convergence is characterized by the relation
$$
\big |(f-R_n)(z)\big |^{1/n}\overset{\mcap }\to e^{-2g_F(z,\infty)}.
$$
\end{theoremSta}

We note again that in the special case when $E=\{\infty\}$,
the external field is absent
(i.e. ${\mathcal V}^{-\delta
_\infty}=0$)
and $f\in\myA (\{\infty\})$, statement $1^\circ$ of Stahl's Theorem yields, in
essence, a positive answer to the existence of
the $\delta_\infty$-minimizing element in the class
$\myk
_{\{\infty\},f}$.

{\bad
It turned out later that the implication $1^\circ\Rightarrow 2^\circ$ in
Stahl's Theorem remains valid in quite more general situations.
E.~Rakhmanov and A.~Martinez-Finkelshtein~\cite{MaRa11} (see also \cite{Rak12},
\cite{MaRaSu12},~\cite{KuSi15}) proved, using the variational method, that the
$\mu$-minimizing element $F$ in $\myk_{E,f}$, which consists of a finite number
of continua is necessarily an $S$-compactum in the external field ${\mathcal V}^{-\mu}$.
In~\cite{MaRa11}, the case of a discrete measure $\mu$ consisting of finitely many
point masses was treated. In the present paper, we consider the case of a
general measure $\mu$. The variational method introduced by Rakhmanov and
Martinez-Finkelshtein is
crucial
for the proof of the following theorem.
}

\begin{theorem}\label{t1}
Let $E=\bigsqcup_{j=1}^mE_j\in {\mathcal E}_m$,
$f\in\myA (E)$, $\mu\in\MM (E)$, $\mu (E_j)>0$, $j=1,\dots ,m. $ Suppose that
the $\mu$-minimizing element $F$ in the class $\myk_{E,f}$ consists of a
finite number of continua.
Then $F$ is symmetric in the external field
${\mathcal V}^{-\mu}$ and
$\overline{\Bbb C}\setminus F=\bigcup_{j=1}^{m}D_j$, where $D_j\supset E_j$. Moreover, the domains
$D_j$ and $D_k$, $j,k=1,\dots ,m$, either do not intersect one another, or coincide.
\end{theorem}

Assume that $K$ is the $\mu$-minimizing
compactum in $\myk_{E,f}$. It is worth noting
that it
is not difficult to find a compactum $F\subseteq K$ consisting of a finite
number of continua and belonging to $\myk_{E,f}$. The particular case when all
compacta $E_j$ are single points $\{ e_j\}$, $j=1,\dots ,m$, was
considered in \cite{Bus15}.

The implication $2^\circ\Rightarrow 3^\circ$ in Stahl's Theorem
can be extended to the more general case
of multipoint Pad{\'e} approximants; for
the definition
see below.

Let $\myt{E}_n=\{ e_{n,1},\dots ,e_{n,2n+1}\}$ be a point set in the extended
complex plane $\overline{\Bbb C}$ (each point preserves its multiplicity),
and let
$$
\omega_{\myt{E}_n}(z)=\prod_{|e_j|\leq 1}(z-e_j)\prod_{|e_j|>1}(1-z/e_j)
$$
be a spherically normalized polynomial associated with the set $\myt{E}_n$. Suppose that $f\in\HH(\myt{E}_n)$.


The multipoint Pad{\'e} approximant of order $n$ of the function $f$
(at the points in $\myt{E}_n$) is the rational function
$R_n=P_{n}/Q_{n}$ determined by ~\eqref{h1} and by
\begin{equation}
\frac{fQ_{n}-P_{n}}{\omega_{\myt{E}_n}}(z)\in\HH(\myt{E}_n\cap\Bbb C),
\label{Q2}
\end{equation}
\begin{equation}
z^{n+1}\frac{fQ_{n}-P_{n}}{\omega_{\myt{E}_n}}(z)\in\HH(\myt{E}_n\cap\{\infty\})
\label{Q4}
\end{equation}
(if $\myt{E}_n\cap\Bbb C=\varnothing$ or $\myt{E}_n\cap\{\infty\}=\varnothing$, then condition~\eqref{Q2} or~\eqref{Q4}, respectively, fails).

The special case $e_{1,n}=\dots=e_{2n+1,n}$ coincides with the concept of classical Pad{\'e} approximant of $f$ at infinity, if $e_{1,n}=\infty$ (or at zero, if
$e_{1,n}=0$).

Everywhere below we
assume that $f\in\HH(E)$. Furthermore, we assume that the sets $\myt{E}_n=\{
e_{n,1},\dots ,e_{n,2n+1}\}$, $n=1,2,\dots$, of the nodes of interpolation are
located on $E$, and their limit distribution is described by
$$
\frac 1{2n}\delta_{\myt{E}_n}\overset{*}\to\mu,\quad n\to\infty.
$$
For tables of interpolation nodes, satisfying
the above
conditions we adopt the notation
$(E,\mu)$.

Let $F\in{\mathcal F}$, and $\Omega$ be a neighborhood of the compactum $F$. We denote
${\HH}_0(\Omega\setminus F)$ the set of all functions $f$
which are
holomorphic in
$\Omega\setminus F$ (or piecewise holomorphic, if
$\Omega\setminus F$ is not connected), and
possess
continuous limit values from both sides
on the arcs of $F_0\setminus B$, where $B=B(f)$ is
a certain
compactum of zero
capacity, wherein the jump $\chi_f$ of $f$ does not have zeros on $F_0\setminus
B$.

A.~A.~Gonchar and E.~A.~Rakhmanov proved~\cite{GoRa87} a fundamental theorem
(see also~\cite{Gon03})
on the
asymptotic behavior
of the zeros of polynomials, satisfying non-Hermitian orthogonality
conditions
on a
compactum
which
possesses the symmetry property.
The notion of symmetry of a compactum $F\in {\mathcal F}$ was introduced in
the present paper only for the case of harmonic fields of the form
$\psi (z)={\mathcal
V}^{-\mu }(z)$ with $S(\mu)\cap F=\varnothing$. For this reason, we formulate
Gonchar--Rakhmanov's theorem only for harmonic fields of
such type.

\begin{theoremGoRa}
Let $F\subset\Bbb C$ be a compactum of positive capacity, $\Omega$ a neighborhood of $F$, and $\mu\in\MM
(\overline{\Bbb C}\setminus\Omega)$. Assume that the sequence of functions $\Psi_n$ and the function $f$
satisfy
the following conditions:

$1^\circ$. $\Psi_n$ are holomorphic in $\Omega, n =1,2,\dots$, and
$$
\psi_n(z)=\frac{1}{2n}\log\frac{1}{|\Psi_n(z)|}\rightrightarrows\psi (z)={\mathcal V}^{-\mu }(z),
\quad z\in\Omega.
$$

$2^\circ$. $F$ is symmetric in the external field ${\mathcal V}^{-\mu }$.

$3^\circ$. The complement $\overline{\Bbb C}\setminus F$ is connected.

$4^\circ$. $f\in {\HH}_0(\Omega\setminus F)$.

If the polynomials $Q_n$, $\deg Q_n\leq n$ ($Q_n\not\equiv 0$) satisfy
the orthogonality conditions
\begin{equation}
\oint_FQ_n(t)\Psi_n(t)f(t)t^\nu dt=0,
\quad\nu =0,1,\dots,n-1\quad (n=1,2,\dots),
\label{F216}
\end{equation}
then the following statements are true:

(i) $\frac 1n\delta_{Q_{n}}\overset{_*}\to\widetilde{\mu }_F$ as $n\to\infty$,
where $\widetilde{\mu }_F$ stands for the balayage of $\mu $ onto $F$.

(ii) if
$Q_n$ are spherically normalized, then
\begin{equation}
\biggl|\oint_F\frac{Q_n^2(t)\Psi_n(t)f(t)dt}{t-z}\biggl |^{1/n}\overset{\mcap}
\to e^{-2w_F^{\mu}},\quad z\in\overline{\Bbb C}\setminus F,
\label{Q668}
\end{equation}
where $w_F^{\mu }$ is the constant from~\eqref{pt77}.
\end{theoremGoRa}

As a corollary, Gonchar and Rakhmanov deduced the following result.

\begin{theoremGoRacon}
Let $F$ be a compactum symmetric in the field ${\mathcal V}^{-\mu}$,
where $\mu\in\MM (E)$ and $E\cap F=\varnothing$. Suppose that the complement $\overline{\Bbb
C}\setminus F$ is connected and let $\{ R_n\}_{n=1}^\infty$ be the sequence of multipoint Pad{\'e} approximants $P_n/Q_n$ of the function $f\in
{\HH}_0(\overline{\Bbb C}\setminus F)$, associated with the $(E,\mu)$-table of nodes of interpolation.Then

(i) $\frac 1n\delta_{Q_{n}}\overset{_*}\to\widetilde{\mu}_F$, where $\widetilde{\mu}_F$ means the balayage of
$\mu$ onto $F$.

(ii) $R_{n}(z)\overset{\mcap}\rightarrow f(z),\quad z\in\overline{\Bbb C}
\setminus F$. The rate of convergence is given by
$$
\big |(f-R_n)(z)\big |^{1/n}\overset{\mcap}\rightarrow e^{-2G_F^\mu (z)}.
$$

\end{theoremGoRacon}

Subsequently, it was proved~\cite{Bus13} that under the condition $\psi(z)={\mathcal V}^{-\mu}(z)$, the
requirement $3^\circ$
in Gonchar--Rakhmanov's theorem
concerning
connectivity of the complement $\overline{\Bbb C}\setminus F$,
can be weakened  to the
condition
that each connected component of the complement $\overline{\Bbb
C}\setminus F$ contains a nonempty inner boundary arc.
For details, the reader is referred to~\cite{Bus13} and~\cite{Bus15}.

Theorem~\ref{t1} and the Corollary of Gonchar--Rakhmanov's theorem are
focused
on the fact that
the positive solution of the $\mu$-minimization problem in the set $\myk_{E,f}$
is of crucial importance
in the theory of rational approximations.
Below
we list some cases
when a positive answer
has been already
achieved.
In each of these cases, the following line of reasoning is used.

Let $F_n$, $n=1,2,\dots$ be a sequence of compacta in
$\myk_{E,f},$ each of them being a union of a finite number (independent on $n$) of continua, and such that
$$
\lim_{n\to\infty}\mcap_{\mu }F_n=\inf_{K\in\myk_{E,f}}\mcap_{\mu }K.
$$
Suppose that there is a mapping $T\colon\Bbb C\to\Bbb C$ with the property
$$
T(F_n)\in\myk_{E,f}\, ,\quad \mcap_{\mu }T(F_n)\leq\mcap_{\mu }F_n\quad \text{and}\quad T(F_n)\subset\Bbb C\setminus\Omega,
$$
where $\Omega$ is some neighborhood of the compactum $E$ which is independent on $n$.
Then, extracting from the sequence
$\{ T(F_n)\}_{n=1}^\infty$
a subsequence which converges to some compactum
$F=\lim_{n\in\Lambda}T(F_n)$ in the Hausdorff metric, and applying
equality~\eqref{h3}, we obtain the equality
$$
\mcap_{\mu }F=\lim_{n\in\Lambda}\mcap_{\mu }T(F_n)=\inf_{K\in\myk_{E,f}}\mcap_{\mu }K.
$$
It is not difficult to check that
$F\in\myk_{E,f}$. Thus, the compact $F$
solves the $\mu$-minimization problem in the set $\myk_{E,f}$.

While proving the first statement of his theorem, H.~Stahl showed that if all
singularities of $f\in\myA (\{\infty\})$ are located in the disk
$U_{R}=\{ |z|\leq R\}$, then as the mapping $T$ one can take the
radial projection $T_R$ onto the disk $U_{R}$; that is
$$\begin{array}{ll}
T_Rz=z,& \hbox{for}\, |z|\leq R,\\
T_Rz=z\frac{R}{|z|},& \hbox{for}\,|z|>R.\\
\end{array}
$$

The case of the so-called two-point Pad\'{e} approximants when $E=(\{ 0\}\cup\{\infty\})\in
{\mathcal E}_2\,$, $f\in\myA (E)$ and the measure $\mu
=\frac{\delta_0+\delta_\infty}2$ was considered in~\cite{BuMaSu12}. It was
shown that as the mapping $T$ one can
take the radial projection onto the annulus $U_{r,R}=\{ r\leq |z|\leq
R\}$ ($0<r< R<\infty$), which contains the set $A_f$ of the singularities of
$f$. The mapping $T=T_{r,R}$ is defined
as follows
$$
\begin{array}{lll}
T_{r,R}z=z, &\hbox{ for} \,z\in U_{r,R}, \\
T_{r,R}z=z\frac{r}{|z|}, &\hbox{ for}\, |z|<r,\\
T_{r,R}z=z\frac{R}{|z|}, &\hbox{ for}\, |z|>R.\\
\end{array}
$$

The proposition used in~\cite{BuMaSu12} on the decay of $\mcap_{\mu }K$ under radial projection onto
the annulus $U_{r,R}$ is
no longer true
if $\mu\neq\frac{\delta_0+\delta_\infty}2.$
Nevertheless, it was shown in~\cite{Bus15} that the posed problem is also solvable, if
$E=\bigsqcup_{j=1}^m\{ e_j\}\in {\mathcal E}_m$, $f\in\myA (E)$,
$\mu =\sum_{j=1}^mp_j\delta_{e_j}\in\MM (E)$ ($p_j\geq0$, $\sum p_j=1$).

The solution of the problem is also positive in the case when $E\in {\mathcal
E}_1$, $E\subseteq\{
|z|\geq 1\}$, $\mu\in\MM (E)$ and $f\in\myA (E)$ such that
$A_f\subset U_{R}=\{ |z|\leq R\}$, $R<1$. More exactly, L.~Baratchart,
H.~Stahl and M.~Yattselev
established~\cite{BaStYa12} that the mapping $T_R$ onto the disk $U_{R}$, as defined above
possesses the desired
properties.
To be precise, we present their result in~\cite{BaStYa12} in the form that
is suitable for the purpose of the current paper.

\begin{theoremBaStYa}
Let
$E\in {\mathcal E}_1,$
$E\subseteq\{ |z|\geq 1\}$, $\mu\in\MM (E)$, $f\in\myA (E)$ and
$A_f\subset U_{R}=\{ |z|\leq R\}$, where $R<1$. Then there is a
$\mu$-minimizing element $F$ in the family $\myk_{E,f}$ and $F\subseteq U_{R}$.
This $\mu$-minimizing
element $F$ is symmetric in the external field $\mathcal V^{-\mu}$, and $F$ is a union
of closures of the critical trajectories of a quadratic differential.
\end{theoremBaStYa}

Note that in a general case
the minimization problem in $\myk_{E,f}$ with respect to the measure
(with no additional assumptions) may not be
solvable.
This may happen in the case of $E\in {\mathcal E}_m$, $\mu\in\MM (E)$ and $f\in\myA (E)$, selected in a special way. In fact,
it is easy to
construct a continuum $E\in {\mathcal E}_1$, a measure $\mu\in\MM (E)$
and a function $f\in\myA (E)$ such that the inequality
\begin{equation}
\mcap_\mu F>\inf_{K\in\myk_{E,f}}\mcap_\mu K.
\label{3}
\end{equation}
holds for every $F\in\myk_{E,f}. $

In the present paper, we show that such an example applies not only to the case
when $E\in {\mathcal E}_1$ and $\mu\in\MM (E)$
selected in a special way, but
also to the concrete simple continuum $E=[a,b]$, for some
sufficiently wide class
of natural measures supported on $[a,b]$, and appropriate functions
$f\in\myA (E)$.

More
precisely
let $E=[a,b]\subset\mathbb C$.
Denote
$\myt{\MM}(E)$
the subclass of all measures $\mu$ from class $\MM (E)$, such that the
logarithmic potential of $\mu$, $V^\mu (z):=-\int\log |z-t|d\mu (t)$, is
continuous at some inner point $z_0=z_0(\mu)$
of the closed interval $E$, and the condition
\begin{equation}
\lim_{t\to 0}\mu (\{|z-z_0|\leq t\})=0
\label{45}
\end{equation}
is satisfied.
We can easily see that both the equilibrium (Chebyshev) measure
for $E$, and the normalized Lebesgue measure,
possess these
properties,
and thus they
are both
contained in $\myt{\MM}(E)$.

\begin{theorem}\label{t2}
Let $E=[a,b]\subset\mathbb C$ be a closed interval, $\mu\in\myt{M}(E)$.
There exits a function $f\in\myA (E)$ such that the $\mu$-minimization
problem~\eqref{q1} in the family of compacta $\myk_{E,f}$ is not solvable.
\end{theorem}

{\bad
To conclude this section we
remark
that in view of Theorem~\ref{t1} the
$\mu$-minimizing compactum from the set $\myk_{E,f}$ with respect to the measure
$\mu$ possesses an $S$-property. The existence of a compactum with the
$S$-property is
crucial
not only
in the study of distribution of zeros
of orthogonal polynomials, but also
in
the derivation of the formulas of strong asymptotics, on the basis of the matrix
Riemann--Hilbert method, and other methods as well
(see~\cite{Nut86},~\cite{Sta96},~\cite{Bus01},~\cite{Apt08},
\cite{GoRaSu11},~\cite{DeaKui11},~\cite{ApYa11},~\cite{MaRaSu12},~\cite{RaSu13},
~\cite{KoSu13}).
}

The results of the present paper
were
partly announced in~\cite{BuSu14}.

\section{Proof of Theorem~\ref{t1}}\label{s4}

Under the conditions of Theorem~\ref{t1}, the $\mu$-minimizing element $F$ in the class
$\myk_{E,f}$, where $E=\bigsqcup
_{j=1}^mE_j\in {\mathcal E}_m$, $\mu\in\MM (E)$, consists of a finite number of continua. Let
$A_f=\bigcup_{j=1}^mA_{f_j}=\{ a_1,\dots ,a_p\}$ be the set of all singularities of $f\in\myA (E)$.

Fix $w\in\Bbb C\setminus F$ and write $h_w(z)=\frac{A_p(z)}{z-w}$,
where $A_p(z)=\prod_{l=1}^p(z-a_l)$.
Applying
standard arguments, we define the mapping (``variation'') $z\mapsto z_t^h=z+th_w(z)$
(see~\cite{PeRa94},~\cite{MaRa11},~\cite{MaRaSu11b}), where $t\in\Bbb C$ is a complex-valued parameter.
If $t$ is small enough, say $0<|t|\leq\varepsilon_0$, then the mapping is univalent in a neighborhood of $F$. Furthermore, $F\mapsto\{
z_t^h,\, z\in F\} = F_t^h$, the set $A_f$ of the singularities of $f$ remains
stable under the mapping
$z\mapsto z_t^h$, and the measure $\nu\in\MM(F)$ maps to the measure
$\nu_t^h\in\MM(F_t^h)$ where $\nu_t^h(B_t^h)=\nu (B)$ for every $B_t^h\subset
F_t^h$. Thus, $d\nu_t^h(z_t^h)=d\nu (z)$. Taking into account the definition
of the class $\myA(E)$, the stability of the set $A_f$ under the mapping
$z\mapsto z_t^h$ and the condition that $F$ consists of a finite number
of continua, we can conclude that $F_t^h\in\myk_{E,f}$.

We remind the definition of the energy $I_\mu (\nu)$ of a measure
$\nu\in\MM(F)$ in the presence of the external field ${\mathcal V}^{-\mu }$
(i.e., the $\mu$-weighted energy):
\begin{equation}
I_\mu (\nu):=
-\iint\log|z-\zeta |d\nu (z)d\nu (\zeta)+2\int {\mathcal V}^{-\mu }(z)d\nu (z).
\label{AB1}
\end{equation}
In accordance with this definition, and
that
of the measure $\nu_t^h$
we have
$$
I_\mu (\nu_t^h)=
-\iint\log|z-\zeta |d\nu_t^h(z)d\nu_t^h(\zeta)+2\int {\mathcal V}^{-\mu }(z)d\nu_t^h(z)
$$
$$
=
-\iint\log|z_t^h-\zeta_t^h|d\nu (z)d\nu (\zeta)+2\int {\mathcal V}^{-\mu }(z_t^h)d\nu (z).
$$
Subtracting equality~\eqref{AB1} from the latter, we get a formula which
describes the increment of the energy of the measure $\nu$ in the field
${\mathcal V}^{-\mu }$ under the mapping $z\mapsto z_t^h$:
\begin{gather}
I_\mu (\nu_t^h)-I_\mu (\nu)=
-\iint\log\biggl|\frac {z_t^h-\zeta_t^h}{z-\zeta }\biggr|d\nu (z)d\nu (\zeta)+
2\iint\log\biggl|\frac{z_t^h-\zeta}{z-\zeta}\biggr|d\mu (\zeta)d\nu (z)\notag\\
=\Re\biggl\{ -\iint\log\biggl (1+\frac {t\big (h_w(z)-h_w(\zeta)\big)}{z-\zeta }\biggr)d\nu (z)d\nu (\zeta)\notag\\
+2\iint\log\biggl(1+\frac{th_w(z)}{z-\zeta }\biggr)d\mu (\zeta)d\nu (z)
\biggr\}=
\Re\biggl\{ -t\iint\frac {\big (h_w(z)-h_w(\zeta)\big)}{z-\zeta }d\nu (z)d\nu (\zeta)\notag\\
+2t\iint\frac{h_w(z)}{z-\zeta }d\mu (\zeta)d\nu (z)\biggl\}+O(t^2)=\Re tH_{w,\mu }(\nu)+O(t^2),
\label{pt161}
\end{gather}
where
\begin{equation}
H_{w,\mu }(\nu)=-\iint\frac{h_w(z)-h_w(\zeta)}{z-\zeta}d\nu (\zeta)d\nu (z)+2\iint\frac{h_w(z)}{z-\zeta }d\mu (\zeta)d\nu (z).
\label{A7}
\end{equation}

Applying equality~\eqref{pt161} to the measures $\widetilde{\mu}_F$ and $\sigma_{t,h}$, where
$\sigma_{t,h},$ in the class $\MM (F)$ is such that $(\sigma_{t,h})_t^h=\widetilde
{\mu}_{F_t^h}$, we obtain two equalities
\begin{equation}
I_\mu\big ((\widetilde {\mu}_{F})_t^h\big)-I_\mu\big (\widetilde {\mu}_{F}\big)=\Re tH_{w,\mu } (\widetilde {\mu}_F)+O(t^2)=O(t),
\label{A1}
\end{equation}
\begin{equation}
I_\mu\big (\widetilde {\mu }_{F_t^h}\big)-I_\mu (\sigma_{t,h})=\Re tH_{w,\mu }(\sigma_{t,h})+O(t^2)=O(t).
\label{A2}
\end{equation}

It is well known~\cite{SaTo97} that for each compactum $K$
whose intersection
with $S(\mu)$
is empty,
the following equalities are valid:
\begin{equation}
I_\mu (\widetilde {\mu}_{K})=\inf_{\nu\in\MM (K)}I_\mu (\nu)
\label{A221}
\end{equation}
and
\begin{equation}
\mcap_\mu K=e^{-I_\mu (\widetilde {\mu}_{K})}.
\label{A22}
\end{equation}
Therefore,
\begin{equation}
I_\mu (\sigma_{t,h})\geq I_\mu (\widetilde {\mu}_{F})\quad\text{and}\quad
I_\mu ((\widetilde {\mu}_{F})_t^h)\geq I_\mu (\widetilde {\mu}_{F_t^h}).
\label{A23}
\end{equation}
Using~\eqref{A2}, the second inequality in~\eqref{A23} and
~\eqref{A1}, we arrive at
$$
I_\mu (\sigma_{t,h})=I_\mu (\widetilde {\mu}_{F_t^h})+O(t)\leq I_\mu ((\widetilde {\mu}_F)_t^h)+O(t)=
I_\mu (\widetilde {\mu}_{F})+O(t),
$$
Combining this inequality and the first inequality in~\eqref{A23} we
obtain, by letting
$t\to 0$
\begin{equation}
I_\mu (\sigma_{t,h})\to I_\mu (\widetilde {\mu}_{F}).
\label{A31}
\end{equation}

It is well known that $\widetilde {\mu}_{K}$ is the unique measure in $\MM
(K)$ for which~\eqref{A221} holds.
Taking into account this observation, the principle of descent~\cite[chapter~I,
\S\,3, Theorem~1.3]{Lan66} and~\eqref{A31}, we deduce that
$\sigma_{t,h}$ converges weakly to the measure $\widetilde
{\mu}_{F}$ as $t\to 0$. This yields
\begin{equation}
\lim_{t\to 0} H_{w,\mu } (\sigma_{t,h})=H_{w,\mu } (\widetilde {\mu}_{F}).
\label{A32}
\end{equation}

As noticed previously, $F_t^h\in\myk_{E,f}$ for all $t$ small enough. Under the conditions of
Theorem~\ref{t1}, the compactum $F$ is the $\mu$-minimizing element in
$\myk_{E,f}$. Thus $\mcap_\mu F\leq\mcap_\mu F_t^h$, which
thanks to~\eqref{A22}, is equivalent to the inequality $I_\mu (\widetilde {\mu}_{F})\geq I_\mu (\widetilde
{\mu}_{F_t^h})$. Using this inequality, the first inequality in~\eqref{A23} and
equality~\eqref{A2}, we obtain
$$
0\geq I_\mu (\widetilde {\mu}_{F_t^h})-I_\mu (\widetilde {\mu}_{F})\geq I_\mu (\widetilde
{\mu}_{F_t^h})-I_\mu (\sigma_{t,h})=\Re tH_{w,\mu }(\sigma_{t,h})+O(t^2),
$$
Since the parameter $t\to 0$ is arbitrary,
the
latter
inequality is possible only if
$\lim_{t\to 0} H_{w,\mu }(\sigma_{t,h})=0$. From here and from~\eqref{A32}
we obtain
\begin{equation}
\label{A33}
H_{w,\mu } (\widetilde {\mu}_{F})=0.
\end{equation}

Applying the definition~\eqref{A7} to $\nu =\widetilde {\mu}_{F}$ on the left hand
side of equality~\eqref{A33}, as well as
the
explicit representation
of the function $h_w(z)=\frac{A_p(z)}{z-w}$,
we rewrite~\eqref{A33} as follows:
\begin{gather}
-\iint\biggl(\frac {A_p(z)}{(z-w)(z-\zeta)}-\frac{A_p(\zeta)}
{(\zeta -w)(z-\zeta)}\biggr)d\widetilde{\mu }_F (z)d\widetilde{\mu }_F (\zeta)+
\notag\\
+2\iint\frac{A_p(z)}{(z-w)(z-\zeta)}d\mu (\zeta)d\widetilde{\mu }_F (z)\equiv0,
\quad w\in\CC\setminus(E\cup F).
\label{t81}
\end{gather}

We note that
$$
A_p(z)-A_p(w)=(z-w)A_{p,1}(z,w),
$$
where $A_{p,1}(z,w)$ is a polynomial of degree $p-1$
in
each variable $z,w$.
Therefore,
\begin{equation}
\frac{A_p(z)}{(z-w)(z-\zeta)}
=\frac{A_p(w)}{(z-w)(z-\zeta)}+\frac{A_{p,1}(z,w)}{z-\zeta }=\frac{A_p(w)}{\zeta -w}\biggl (\frac 1{z-\zeta }-\frac 1{z-w}\biggl)+\frac{A_{p,1}(z,w)}{z-\zeta }.
\label{t84}
\end{equation}

We also note that the expression
$$
A_p(z)(w-\zeta)+A_p(w)(\zeta -z)+A_p(\zeta)(z-w)
$$
is a polynomial of degree $p$
in
each variable $z,w,\zeta$,
with zeros
at $z=w$, $z=\zeta$, $w=\zeta$. Thus
$$
\frac{A_p(z)}{(z-w)(\zeta -z)}+\frac{A_p(w)}{(w-\zeta)(z-w)}+\frac{A_p(\zeta)}{(\zeta -z)(w-\zeta)}=A_{p,2}(z,w,\zeta),
$$
where $A_{p,2}(z,w,\zeta)$ is a polynomial of degree $p-2$
in
each $z,w,\zeta$. Consequently,
\begin{equation}
\frac {A_p(z)}{(z-w)(z-\zeta)}-\frac{A_p(\zeta)}{(\zeta -w)(z-\zeta)}=-A_{p,2}(z,w,\zeta)-\frac{A_p(w)}{(z-w)(\zeta -w)}.
\label{t83}
\end{equation}

For $w\in\myo\CC\setminus(E\cup F)$ we obtain,
by substitution of
~\eqref{t84} and~\eqref{t83} into~\eqref{t81},
$$
0=\iint\biggl(A_{p,2}(w,z,\zeta)+\frac{A_p(w)}{(z-w)(\zeta -w)}\biggr)d\widetilde{\mu }_F (z)d\widetilde{\mu }_F (\zeta)+
$$
$$
2\iint\biggl(\frac{A_p(w)}{\zeta -w}\biggl(\frac 1{z-\zeta }-\frac 1{z-w}\biggr)+\frac{A_{p,1}(w,z)}{z-\zeta }\biggr)d\mu (\zeta)d\widetilde{\mu }_F (z)
$$
$$
= A_{p,3}(w)+A_p(w)\varphi_{\widetilde{\mu }_F}(w)^2+2A_p(w)\int\frac{\varphi_{\widetilde{\mu }_F}(\zeta)}{\zeta -w}d\mu (\zeta)-2A_p(w)\varphi_{\widetilde{\mu }_F}(w)\varphi_{\mu}(w)
,
$$
where
$$
A_{p,3}(w)=\iint A_{p,2}(w,z,\zeta)d\widetilde{\mu }_F (z)d\widetilde{\mu }_F (\zeta)+2\iint\frac{A_{p,1}(w,z)}{z-\zeta }d\mu (\zeta)d\widetilde{\mu }_F (z)
$$
is a polynomial of degree not greater
than
$p-1$. The integrals
$$
\varphi_{\widetilde{\mu }_F}(w)=\int\frac {d\widetilde{\mu }_F (z)}{z-w}
\quad\text{and }\quad\varphi_{\mu}(w)=\int\frac {d\mu (z)}{z-w},\quad
w\in\myo\CC\setminus(E\cup F).
$$
represent the Cauchy transformations of the measures $\widetilde{\mu
}_F $ and $\mu$, respectively.

Rewrite the equality which we obtained as follows:
\begin{equation}
\bigl(\varphi_{\widetilde{\mu }_F}(w)
-\varphi_{\mu}(w)\bigr)^2=\varphi_{\mu}(w)^2-\frac{A_{p,3}(w)}{A_p(w)}-2\int\frac{\varphi_{\widetilde{\mu }_F}(\zeta)}{\zeta -w}d\mu (\zeta).
\label{t91}
\end{equation}
Denote the right hand side of~\eqref{t91} by $R(w).$ The function
$R(w)$ is meromorphic on $F$, and has poles only at the zeros of the polynomial $A_p(w)$.

By virtue of~\eqref{t91},
we derive for the multi-valued complex potential $\myVV^{\mu -\widetilde{\mu }_F} (z)$ of the charge $\mu
-\widetilde{\mu }_F$
$$
\myVV^{\mu -\widetilde{\mu }_F}(z):=\int\log(z-\zeta)d(\widetilde{\mu }_F -\mu)(\zeta)=
$$
$$
\int\biggl(\int ^z\frac 1{w-\zeta }dw\biggl)d(\widetilde{\mu }_F -\mu)(\zeta)=\int ^z\big (\varphi_{\mu}(w)-\varphi_{\widetilde{\mu }_F}(w)\big)dw=\int ^z\sqrt{R(w)}dw.
$$
Hence, for the corresponding real parts we have
\begin{equation}
{\mathcal V}^{\mu -\widetilde{\mu }_F}(z)=\Re\int ^z\sqrt{R(w)}dw\enspace .
\label{t93}
\end{equation}

It follows from~\eqref{pt77} that ${\mathcal V}^{\mu -\widetilde{\mu }_F}(z)$ is,
up to a constant $w_F^\mu$, identical to $G_F^\mu (z)$.

Hence, by~\eqref{t93} and
in view
of the fact that $G_F^\mu =0$ on $F$,
$$
G_F^\mu (z)=\Re\int_a^z\sqrt{R(w)}dw\, ,\quad\text{where}\quad a\in F\enspace .
$$
Since $F=\overline{\{ z:\, G_F^\mu (z)=0\}}$,
we have
$$
F=\overline{\biggl\{z:\,\Re\int_a^z\sqrt{R(w)}dw=0\biggr\}}.
$$
Thus, $F$ consists of a finite number of analytic arcs, the endpoints of which
belong to the union of the set $\{a_1,\dots ,a_p\}$ and the set of zeros
of the function $R$.

Below
we prove that the compactum $F$ is symmetric in the external field ${\mathcal V}^{-\mu}$
(or, in
the
other words, we establish equality~\eqref{S2}). Let $\zeta\in F_0$. Denote
$\nabla_\pm
{\mathcal V}^{\mu -\widetilde{\mu }_F}$ the limit values of the gradient of the potential ${\mathcal V}^{\mu -\widetilde{\mu }_F}$.
Since ${\mathcal V}^{\mu
-\widetilde{\mu }_F}\equiv\const$ for $z\in F$, then $F$ turns out to be a level curve
of ${\mathcal V}^{\mu-\widetilde{\mu }_F} $ (i.e., an equipotential curve),
which
implies
$$
\biggl|\frac{\partial {\mathcal V}^{\mu -\widetilde{\mu }_F}}{\partial n_\pm
}(\zeta)\biggr|=\biggl|\nabla_{\pm}{\mathcal V}^{\mu -\widetilde{\mu }_F}(\zeta)\biggr|,\quad\zeta\in F_0.
$$
On the other hand, ${\mathcal V}^{\mu -\widetilde{\mu }_F} =\Re\myVV^{\mu -\widetilde{\mu
}_F}$ and $\dfrac {d\myVV^{\mu -\widetilde{\mu }_F}}{d\zeta}(\zeta)
=\sqrt{R(\zeta)}$. Hence,
$$
\big |\nabla {\mathcal V}^{\mu -\widetilde{\mu }_F}(\zeta)\big |=
\biggl |\frac {d\myVV^{\mu -\widetilde{\mu }_F}}{d\zeta}(\zeta)\biggl |=|\sqrt{R(\zeta)}|,\quad\zeta\in\Bbb C\setminus F.
$$
Therefore,
$\bigl|\nabla_+{\mathcal V}^{\mu -\widetilde{\mu }_F}(\zeta)\bigr|=\big
|\nabla_-{\mathcal V}^{\mu -\widetilde{\mu }_F}(\zeta)\big |$ which in turn yields that
$$
\biggl|\frac{\partial {\mathcal V}^{\mu -\widetilde{\mu }_F}}{\partial n_+}(\zeta)\biggl |=\biggl |\frac{\partial {\mathcal V}^{\mu -\widetilde{\mu }_F}}{\partial n_-}(\zeta)\biggr|,\quad\zeta\in F_0.
$$
Since
$G_F^\mu (z)=0$ for $z\in F$ and $G_F^\mu (z)>0$ for
$z\in\Bbb C\setminus F$, we have
$\dfrac{\partial{\mathcal V}^{\mu -\widetilde{\mu}_F}}{\partial n_\pm}(\zeta)>0$,
$\zeta\in F_0$. From here we derive the
equality
$$
\frac{\partial {\mathcal V}^{\mu -\widetilde{\mu }_F}}{\partial n_+}(\zeta)=\frac{\partial {\mathcal V}^{\mu -\widetilde{\mu }_F}}{\partial n_-}(\zeta),\quad\zeta\in F_0,
$$
which coincides with~\eqref{S2}.

Now, let the connected component $G$ of the complement of $F$ be such that
$G\cap E=G\cap S(\mu)=\varnothing$. Then the function on the left hand side of
the equality~\eqref{t93} is harmonic in the domain $G$ and equals a constant on
$\partial G\subset F$. Thus this function is identically  constant on $G$. The
last statement contradicts the right hand side of~\eqref{t93}. Consequently,
each connected component of $F$ has necessarily a nonempty intersection with
$E$. Since $E=\bigsqcup
_{j=1}^mE_j\in {\mathcal E}_m$, then each connected component contains entirely one or several continua $E_j$, $j=1,\dots ,m$.

The proof of Theorem~\ref{t1} is complete.

\section{Proof of Theorem~\ref{t2}}\label{s3}

Without loss of generality, we may assume
that $E=[a,b]$ is located on the real
axis, and $a\leq -1$, $b=1$, $z_0=0$.
Given $k\in\Bbb N$, the linear function $z\mapsto kz$ maps $E$ onto the closed
interval
$E_k=[ka,k]$, where $ka\leq -k$.

By this,
the measure $\mu\in\myt\MM (E)$
transfers
to a measure
$\mu_k\in\MM (E_k)$, where $\mu_k(B)=\mu (\{z: kz\in
B\})$;
therefore,
$d\mu_k(kz)=d\mu (z)$.

Consider
the multi-valued analytic function
$$
f(z)=\sqrt{(z-a_1)(z-a_2)(z-a_3)^{-1}(z-a_4)^{-1}},
$$
where
$$
a_1=\frac{-2+3i}{16},\quad a_2=\frac{2+3i}{16},\quad a_3=\frac{-2-i}{16},
\quad a_4=\frac{2-i}{16}.
$$
Since for the branch points of $f$ we have $\Im a_1,\Im a_2>0$
and $\Im a_3,\Im a_4<0$, we can take a holomorphic branch $f^*$ of $f$ on $\Bbb
R$, determined by the condition $f^*(\infty)=1$. Clearly, $f^*\in\myA (E_k)$
for every $k=1,2,\dots$.

Denote
$L$ the union of the diagonals $[a_1,a_4]$ and $[a_2,a_3]$
of the square with vertices at the points $a_1$, $a_2$, $a_3$, $a_4$. We keep
the notation $f_L$ for the holomorphic branch of $f$
in the domain
$\overline{\Bbb C}\setminus L$
determined by the condition $f_L(\infty)=1$.

We easily see that the proof of Theorem~\ref{t2} will be completed, if we show that there exists a number
$k\in\Bbb N$ such that for all $F\in\myk_{E_k,f^*}$ the inequality
\begin{equation}
\mcap_{\mu_k}F>\inf_{K\in\myk_{E_k,f^*}}\mcap_{\mu_k}K
\label{4}
\end{equation}
holds.
Furthermore, we show that inequality~\eqref{4} is true for all sufficiently
large $k$ and all $F\in\myk_{E_k,f^*}$.

Let
$$
\psi_k(z):=V^{-\mu}(k^{-1}z)-V^{-\mu}(0).
$$
Because of the conditions imposed on the measure $\mu$,
the potential $V^{\mu}(z)$ is continuous at the point $z=0$.
This implies
the uniform convergence of the
sequence of
functions $\psi_k$ on each compact set $K\subset\Bbb C$,
in particular,
on the disk $D=\{ |z|\leq
1\}$, as well:
\begin{equation}
\psi_k(z)\rightrightarrows 0,\quad z\in D,\quad k\to\infty.
\label{144}
\end{equation}
We observe that the spherically normalized potential ${\mathcal V}^{\mu_k}(z)$, given by the equality~\eqref{pt173},
differs from the standard potential $V^{\mu_k}(z)$
by a constant. From here and from $d\mu_k(k\tau)=d\mu (\tau)$, we obtain
$$
{\mathcal V}^{-\mu_k}(z)
=V^{-\mu_k}(z)+C_k=\int\log |z-t|d\mu_k(t)+C_k
$$
$$
=\int\log |k(k^{-1}z-\tau)|d\mu (\tau)+C_k
=\log k+C_k+V^{-\mu}(0)+\psi_k(z)\, .
$$
Now, considering~\eqref{pt172}, we conclude that the inequality~\eqref{4} is equivalent to the inequality
\begin{equation}
\mcap_{\psi_k}(F)>\inf_{K\in\myk_{E_k,f^*}}\mcap_{\psi_k}(K).
\label{44}
\end{equation}
Denoting $\rho_k=\exp\{\max_{z\in D}|\psi_k(z)|\}$, we get for $\rho_k\geq 1$ the estimate
$$
\rho_k^{-1}\leq e^{\psi_k(z)}\leq\rho_k\, ,\, z\in D\,.
$$
Furthermore, by~\eqref{144}, $\lim_{k\to\infty}\rho_k=1$.
From here and from~\eqref{pt172}, after assuming that $K\subseteq D$, we deduce
the
estimates of the capacity $\mcap_{\psi_k} K$
(from below and above) in terms of the standard capacity $\mcap K$, namely

\begin{equation}
\rho_k^{-2}\mcap K\leq\mcap_{\psi_k}(K)\leq\rho_k^{2}\mcap K.
\label{6}
\end{equation}

Assume
that the assertion that the inequality~\eqref{44} holds for all $k$ starting
from
a number $k_0$ and for all $F\in\myk_{E_k,f^*}$, is false.
Then there is an increasing sequence $\{ k_j\}_{j=1}^\infty$ of natural numbers and compacta $F_{k_j}\in\myk_{E_{k_j},f^*}$
such that for all $j=1,2,\dots$
\begin{equation}
\mcap_{\psi_{k_j}}(F_{k_j})
=\inf_{K\in\myk_{E_{k_j},f^*}}\mcap_{\psi_{k_j}}(K).
\label{5}
\end{equation}
Without loss of generality, we may suppose that each compactum $F_{k_j}$
contains at most two continua, and each of these continua
contains a pair of points from the set $\{ a_1,a_2,a_3,a_4\}$,
see~\cite[Lemma~15]{BaStYa12}.

We first establish the inclusion $F_{k_j}\subset D$ for all sufficiently
large $j$. For this purpose, we observe that
$K^* =[a_1,a_2]\cup [a_3,a_4]\in\myk
_{E_{k_j},f^*}$ and that $K^*$ is contained in the disk of radius $2^{-3}\sqrt {2}$ centered at $a_0=\sum_{l=1}^4a_l/4=2^{-4}i$.
Therefore, after taking into account the second inequality in \eqref{6}, and the
equality $\lim_{k\to\infty}\rho_k=1$, for $j$ large enough, say $j\geq j_0$, we derive the estimate
\begin{equation}
\inf_{K\in\myk_{E_{k_j},f}}\mcap_{\psi_{k_j}}(K)\leq\mcap_{\psi
_{k_j}}(K^*)\leq\rho_{k_j}^2\mcap K^*\leq\rho_{k_j}^22^{-3}\sqrt {2}<0.18.
\label{7}
\end{equation}
On the other hand, if the compactum $F_{k_j}\in\myk_{E_{k_j},f^*}$ (consisting of at most
two continua, each of
the latter
containing
two points from
the set
$\{ a_1,a_2,a_3,a_4\}$), is not lying in the disk $D$, then $F_{k_j}$
contains some continuum $F_{k_j}^*$ located in $ D$ and combining some point on the circle $\{ |z|=1\}=\partial D$ with some of the points $a_1$, $a_2$, $a_3$, $a_4$.
The distance from each of these points to $\partial D$ is not smaller than
$1-2^{-4}\sqrt{13}$. Therefore,
$$
\mcap F_{k_j}\geq\mcap F_{k_j}^*\geq 2^{-2}(1-2^{-4}\sqrt{13})> 0.19.
$$
Hence, by the first inequality in~\eqref{6}, for sufficiently large $j$,
say $j\geq j_1$, we have
$$
\mcap_{\psi_{k_j}}(F_{k_j})\geq\mcap F_{k_j}\rho_{k_j}^{-2}> 0.19\rho_{k_j}^{-2}> 0.18.
$$
The
latter
inequality and~\eqref{7} yield for every $j\geq \max\{j_0,j_1\}$
$$
\mcap_{\psi_{k_j}}(F_{k_j})>\inf_{K\in\myk_{E_{k_j},f^*}}\mcap_{\psi_{k_j}}(K),
$$
which
contradicts~\eqref{5}. Consequently,
$F_{k_j}\subset D$ for all $j\geq\max\{j_0,j_1\}$.

From the inclusion $F_{k_j}\subset D$ and from the relation $F_{k_j}\cap
E_{k_j}=\varnothing$ we
obtain,
in particular, that for $j\geq\max\{j_0,j_1\}$ the set $F_{k_j}$ consists exactly of two continua $F_{k_j}^+$ and
$F_{k_j}^-$; the former lies in $D\cap\{\Im z>0\}$ and connects $a_1$ with
$a_2$, the latter is located in
$D\cap\{\Im z<0\}$ and connects $a_3$ with $a_4$.

Recall
the definition of the distance between two compacta $K_1,K_2\subset\overline{\Bbb C}$ in the spherical Hausdorff metric, that is:
\begin{equation}
\rho_H(K_1,K_2)=\inf\{\delta :\, K_1\subset K_2^\delta\quad \text{and}\quad
K_2\subset K_1^\delta\},
\label{K1}
\end{equation}
where
\begin{equation}
K^\delta =\{ z:\, d(z,K)<\delta\},\qquad
d(z,w)=|z-w|(1+|z|^2)^{-1/2}(1+|w|^2)^{-1/2}.
\label{K}
\end{equation}

The notation $K_n\rightarrow K$ stands for convergence of compacta $K_n$,
$n=1,2,\dots$, to the compact $K$ in the spherical Hausdorff metric. We remind that from each sequence of compacta
$K_n\subset\overline{\Bbb
C}$ one can extract a convergent subsequence.

We may suppose without
loss of
generality that the sequence $\{ k_j\}_{j=1}^\infty$
of natural numbers is such that the sequences of continua $\{
F_{k_j}^+\}_{j=1}^\infty$ and $\{ F_{k_j}^-\}_{j=1}^\infty$ converge as
$j\to\infty$ to some continua $F^+$ and $F^-$, respectively (if needed, we pass to subsequences). It can be easily seen that $F:=(F^+\cup F^-)\in\myk_{a_3,a_4}^{a_1,a_2}$, where $\myk
_{a_3,a_4}^{a_1,a_2}$ denotes that set of compacta $K$, consisting of the
union of continua $K^+\cup K^-$, where $K^+\subset D\cap\{\Im
z\geq 0\}$ and connects the point $a_1$ with the point $a_2$, whereas
$K^{-}\subset D\cap\{\Im z\leq 0\}$ and connects $a_3$ with
$a_4$. We show that
\begin{equation}
\mcap F=\inf_{K\in\myk_{a_3,a_4}^{a_1,a_2}}\mcap K.
\label{8}
\end{equation}
In fact, let
$$
S_j^+:\,\{\Im z\geq 0\}\mapsto\{\Im z\geq j^{-1}\}\quad \text{and}\quad
S_j^-:\,\{\Im z\leq 0\}\mapsto\{\Im z\leq -j^{-1}\}
$$
denote the transformation, determined in the following way:
$$
S_j^+z=\Re z+i
\begin{cases}\Im z,&\text{if}\quad\Im z>j^{-1},\\
j^{-1},&\text{if}\quad 0\leq\Im z\leq j^{-1},
\end{cases}
\quad S_j^-z=\myo{S_j^+(\myo z)}.
$$

If $K=(K^+\cup K^-)\in\myk_{a_3,a_4}^{a_1,a_2}$, then we
let
$S_jK=(S_j^+K^+)\cup (S_j^-K^-)$. We easily see that if
$K\in\myk_{a_3,a_4}^{a_1,a_2}$, then $S_jK\in\myk_{E_{k_j},f^*}$ and
$S_jK\to K$ as $j\to\infty$. Hence, from~\eqref{6},
in view of~\eqref{h3} and
that $\lim_{j\to\infty}\rho_{k_j}=1$, we
obtain
that for each
$K\in\myk_{a_3,a_4}^{a_1,a_2}$
$$
\mcap K=\lim_{j\to\infty}\mcap S_jK
=\lim_{j\to\infty}
\mcap_{\psi_{k_j}}(S_jK)\geq\lim_{j\to\infty}\mcap_{\psi_{k_j}}(F_{k_j})
=\lim_{j\to\infty}\mcap F_{k_j}
=\mcap F,
$$
From here, we derive the equality~\eqref{8}.

For $p=4$ and $p=5$ denote by $L_p$ the intersection of $L$ with the halfplane $\{\Im z\geq 2^{-p}\}$.
Recall that $a_0=\sum_{l=1}^4a_l/4=2^{-4}i$.
Therefore, $L_{4}$ is a union of the upper halfdiagonals
$[a_1,a_0]$ and $[a_0,a_2]$ of the square with vertices at the points $a_1,a_2,a_3,a_4$.
We show that $L_{5}\setminus F^+\neq\varnothing$. In fact, otherwise $F^+\supseteq L_{5}\supsetneq L_4$, whence,
$$
\mcap (L_4\cup F^-)<\mcap (L_5\cup F^-)\leq\mcap (F^+\cup F^-)=\mcap F.
$$
This inequality contradicts~\eqref{8}, thanks to the fact that
$(L_4\cup F^-)\in\myk_{a_3,a_4}^{a_1,a_2}$.

Thus, there is a point $z^*\in L_{5}\setminus F^+$. Denote
$\widetilde{\mu
}_{F_{k_j}}$ the balayage of the measure $\mu_{k_j}$
onto $F_{k_j}$. Without
loss of generality,
we may suppose that the sequence of measures $\{\widetilde{\mu
}_{F_{k_j}}\}_{j=1}^\infty$ (passing to subsequences, if needed) converges in the weak sense
\begin{equation}
\widetilde{\mu }_{F_{k_j}}\overset{*}{\rightarrow}\mu^{*}_F,\quad j\to\infty.
\label{18}
\end{equation}
Since $F_{k_j}\to F$, we have $\supp \mu^{*}_F\subset F$. From $z^*\in
L_{5}\setminus F^+\subset L\setminus F$, it follows that
$\mu^{*}_F\neq\lambda^{\equ}_L$, where $\lambda^{\equ}_L$ is the equilibrium measure of the
compactum $L$
(in the absence of an external field).

Fix some $(E,\mu)$-table of points of interpolation $\myt{E}_n=\{
e_{2n+1,l}\}_{l=1}^{2n+1}$, $n=1,2,\dots$.
Then $k_j\myt{E}_n=\{k_je_{2n+1,l}\}_{l=1}^{2n+1}$, $n=1,2,\dots$,
is an $(E_{k_j},\mu_{k_j})$-table.

We get from our assumption~\eqref{5}, after applying Theorem~\ref{t1}
and the Corollary of
Gonchar--Rakhmanov's theorem, that for all $j=1,2,\dots$
\begin{equation}
\frac 1n\delta_{Q_n^{k_j}}\overset{*}{\rightarrow}\widetilde{\mu }_{F_{k_j}},
\quad n\to\infty,
\label{1}
\end{equation}
where $Q_n^{{k_j}}$ is the denominator of the Pad{\'e} approximants of the function $f^*\in\myA(E_{k_j})$, constructed
with respect to the $(E_{k_j},\mu_{k_j})$-table of the interpolation nodes
$k_j\myt{E}_n$.


Let $C(D)$ be the space of all continuous functions on $D$.
Select some countable and everywhere dense set of functions
$g_1,g_2,\dots$ in $C(D)$ ; write, for $\nu_1,\nu_2\in\MM (D)$
$$
\rho (\nu_1,\nu_2):=\sum_{p=1}^\infty\frac 1{2^p\max_{z\in D}|g_p(z)|}\biggl |\int g_pd\nu_1-\int g_pd\nu_2\biggl |.
$$
We note that
\begin{equation}
\rho (\nu_1,\nu_3)\leq\rho (\nu_1,\nu_2)+\rho (\nu_2,\nu_3).
\label{161}
\end{equation}
Moreover, in view of~\eqref{h4} (the definition of weak convergence),
for $\nu_j,\nu\in\MM (D)$
\begin{equation}
\nu_j\overset{*}{\rightarrow}\nu,
\quad j\to\infty\Longleftrightarrow\lim_{j\to\infty}\rho(\nu_j,\nu)=0.
\label{162}
\end{equation}
In particular, relations~\eqref{18} and~\eqref{1} can be rewritten as equalities; namely,
\begin{equation}
\lim_{j\to\infty}\rho (\widetilde{\mu}_{F_{k_j}},\mu^{*}_F)=0\quad
\text{and}\quad\lim_{n\to\infty}
\rho\biggl(\frac 1n\delta_{Q_n^{k_j}},\widetilde{\mu }_{F_{k_j}}\biggr)=0.
\label{163}
\end{equation}
Denote
$J_{n,j}$ the set all
indices
$l\in\{1,\dots ,2n+1\}$ such that
$|e_{2n+1,l}|> k_j^{-1/2}$, and let $|J_{n,j}|$ be it's cardinality.
We get from the definition of $(E,\mu)$-table that for all $n\geq n(j)$
\begin{equation}
|J_{n,j}|/(2n)\geq\mu (E\setminus [-k_j^{-1/2},k_j^{-1/2}])-k^{-1}_j.
\label{28}
\end{equation}

Along with the sequence $\{ k_j\}_{j=1}^\infty$ that we already constructed,
we
introduce a sequence $\{ n_j\}_{j=1}^\infty$ of natural numbers
according to
the
following scheme. We fix
arbitrarily
a starting number $n_1\in\Bbb N$.
Suppose that the numbers $n_1,\dots ,n_{j-1}$
have already been
found, and select the
number $n_j$ in such a way that $n_j>\max\{n_{j-1},n(j)\}$ and
\begin{equation}
\rho\biggl(\frac 1{n_j}\delta_{Q_{n_j}^{k_j}},\widetilde{\mu }_{F_{k_j}}\biggr)
\leq k_j^{-1}.
\label{2}
\end{equation}
Such a choice of the number $n_j$ is possible thanks to the second equality
in~\eqref{163}.

Hence, the numbers $k_j$ and $n_j$ are defined for all $j=1,2,\dots$.
Set $n_j^*=|J_{n_j,j}|$; taking into account the obvious
inequality~$|J_{n_j,j}|\leq 2n_j+1$, we obtain that inequality~\eqref{28} and
condition~\eqref{45} lead to
\begin{equation}
\lim_{j\to\infty} n_j^*/n_j=2.
\label{17}
\end{equation}
Let
$$
T_{n_j}(z)=\prod_{l=1}^{2n_j+1}\biggl (1-\frac{z}{k_je_{2n_j+1,l}}\biggl),
\quad
\omega_{n_j}(z)=\prod_{l\in J_{n_j,j}}\biggl(1-\frac{z}{k_je_{2n_j+1,l}}\biggr).
$$
We note that $\deg\omega_{n_j}=n_j^*$, the polynomial $\omega_{n_j}$ divides $T_{n_j}$, is zero-free for $|z|\leq k_j^{1/2}$
and for each $z\in D$, the inequalities
\begin{equation}
\biggl (1-k_j^{-1/2}\biggl)^{n_j^*}\leq |\omega_{n_j}(z)|\leq\biggl (1+k_j^{-1/2}\biggl)^{n_j^*}
\label{171}
\end{equation}
hold. Thus, the functions
$\Psi_{n_j}(z)=1/{\omega_{n_j}(z)}$, $j=1,2,\dots$, are holomorphic in the disk $D$ and
\begin{equation}
\frac 1{n_j}\log |\Psi_{n_j}(z)|\rightrightarrows 0={\mathcal V}^{-\delta_\infty}(z),
\quad z\in D.
\label{175}
\end{equation}
We remark that the compactum $L$ is symmetric in the external field
${\mathcal V}^{-\delta_\infty}\equiv 0$ (in other words,
in the absence of an external field),
the complement $\overline{\Bbb C}\setminus L$ is connected and $f_L\in
{\HH}_0(\overline{\Bbb C}\setminus L)$. From
the above consideration
it follows that the compact set $L$,
the neighborhood $\Omega =\{ |z|<1\}$ of $L$, the measure $\mu =\delta
_\infty$, the sequence of functions $\Psi_{n_j}$, and the function $f_L$
satisfy conditions $1^\circ$--$4^\circ$ in Gonchar--Rakhmanov's theorem
(after replacing the entire set of natural numbers by the sequence $\{ n_j\}_{j=1}^\infty$).
We later return to this version of Gonchar--Rakhmanov's theorem, a bit stronger
than the original. Before
applying this statement
we
check
whether the polynomials
$Q_{n_j}^{k_j}$ of degree not greater
than
$n_j$
satisfy the following orthogonality conditions:
\begin{equation}
\oint_LQ_{n_j}^{k_j}(t)\Psi_{n_j}(t)f_L(t)t^\nu dt=0,
\quad\nu=0,1,\dots,n_j^*-n_j-2\quad (j=1,2,\dots);
\label{F217}
\end{equation}
where we integrate along a contour encircling $L$ and close enough to it.

In fact, $f_L(z)=f^*(z)=f_{F_{k_j}}(z)$ for all $z\in\myo\CC\setminus D$,
$j\geq\max\{j_0,j_1\}$, where $f_{F_{k_j}}$ is the branch of~$f$ holomorphic in the domain
$\myo\CC\setminus F_{k_j}$ and determined by the condition
$f_{F_{k_j}}(\infty)=1$.
By definition of Pad{\'e} approximants of the function $f_{F_{k_j}}\in
\myA (E_{k_j})$, the difference $\big (Q_{n_j}^{k_j}f_{F_{k_j}}-P_{n_j}^{k_j}\big
)(z)$ vanishes at all zeros of the polynomial $T_{n_j}(z)$ and,
in particular, at all zeros of the polynomial
$\omega_{n_j}(z)$. Because of the fact that all zeros of $\omega_{n_j}(z)$ belong
to the set $\{ |z|>k_j^{1/2}\}\subset\overline{\Bbb C}\setminus D$, on which
$f_{F_{k_j}}=f_L$, we come to the conclusion that
\begin{equation}
\frac{Q_{n_j}^{k_j}f_L-P_{n_j}^{k_j}}{\omega_{n_j}}\in\HH(\overline{\Bbb C}\setminus L).
\label{217}
\end{equation}
The left hand side of the inclusion~\eqref{217} is of order
$O(1/z^{n_j^*-n_j})$ as $z\to\infty$. Therefore, for all $\nu =0,1,\dots ,n_j^*-n_j-2$
\begin{equation}
\oint_{L}\frac{Q_{n_j}^{k_j}f_L-P_{n_j}^{k_j}}{\omega_{n_j}}(t)t^{\nu}dt=0.
\label{a2}
\end{equation}
Since all zeros of the polynomial $\omega_{n_j}(z)$ are lying outside of $D\supset L$, we get
\begin{equation}
\oint_{L}\frac{P_{n_j}^{k_j}}{\omega_{n_j}}(t)t^{\nu}dt=0,\quad\nu =0,1,\dots.
\label{a3}
\end{equation}
Hence, equalities~\eqref{a2} can be rewritten as follows:
$$
\oint_{L}\frac{Q_{n_j}^{k_j}f_L}{\omega_{n_j}}(t)t^{\nu}dt=0,\quad\nu
=0,1,\dots,n_j^*-n_j-2,
$$
which coincide with~\eqref{F217}.

Now, let us address to the relations~\eqref{F217} and~\eqref{F216}.
Notice that in~\eqref{F217} the indices $\nu =n_j^*-n_j-1,\dots ,n_j-1$
are absent. It follows from~\eqref{17} that the number $2n_j-n_j^*+1$
of absent relations does not exceed $o(n_j)$. Taking this into account
and repeating verbatim the proof elaborated by Gonchar and Rakhmanov,
we see that statement $(i)$ of Gonchar--Rakhmanov's theorem follows from
relation~\eqref{F217} just as from relation~\eqref{F216}. Consequently,
by the stronger version of statement $(i)$ in Gonchar--Rakhmanov's
theorem,
\begin{equation}
\frac 1{n_j}\delta_{Q_{n_j}^{k_j}}\overset{*}\to\lambda^{\equ}_L,\quad j\to\infty,
\label{a5}
\end{equation}
where $\lambda^{\equ}_L$ stands for the equilibrium measure of the compactum $L$ (in
the absence of an external field).

On the other hand, inequality~\eqref{161} yields that
$$
\rho\biggl(\frac 1{n_j}\delta_{Q_{n_j}^{k_j}},\mu^{*}_F\biggr)
\leq\rho\biggl(\frac 1{n_j}\delta_{Q_{n_j}^{k_j}},\widetilde{\mu }_{F_{k_j}}\biggr)
+\rho (\widetilde{\mu }_{F_{k_j}},\mu^{*}_F).
$$
Now, considering inequality~\eqref{2} and the first equality in~\eqref{163}, we arrive at the equality
$$
\lim_{j\to\infty}\rho\biggl(\frac 1{n_j}\delta_{Q_{n_j}^{k_j}},
\mu^{*}_F\biggr)=0,
$$
which is equivalent to the relation
\begin{equation}
\frac 1{n_j}\delta_{Q_{n_j}^{k_j}}\overset{*}\to\mu^{*}_F,\quad j\to\infty.
\label{a55}
\end{equation}
As noticed before, $\lambda^{\equ}_L\neq\mu^{*}_F$. Thus,~\eqref{a55}
contradicts~\eqref{a5}. Therefore, the proof of Theorem~\ref{t2} is
complete.

\subsection*{Acknowledgments}
The authors are grateful to the referees for useful suggestions about the manuscript.
This work was supported by the Russian Science Foundation (RSF) under the grant
14-21-00025.












\end{document}